\input ./style/arxiv-ba.cfg
\documentclass[ba,linksfromyear,preprint]{imsart}
\makeatletter
   \@ifpackageloaded{natbib}{}{\usepackage{natbib}}
\makeatother

\pubyear{2015}
\volume{10}
\issue{3}
\firstpage{749}
\lastpage{752}
\doi{10.1214/15-BA944C}% Updated by VTEXPTS2LaTeX.exe, 28.04.2015 08:20
\newcommand{\enquote}[1]{``#1''}

\begin{document}

\begin{frontmatter}
\title{Comment on Article by Ferreira and Gamerman\thanksref{T1}}
\runtitle{Comment on Article by Ferreira and Gamerman}

% Comment, Discussion, Rejoinder
\relateddois{T1}{Main article DOI: \relateddoi[ms=BA944]{Related item:}{10.1214/15-BA944}.}

\begin{aug}
\author[addr1]{\fnms{James V.} \snm{Zidek}\corref{}\ead[label=e1]{jim@stat.ubc.ca}}

\runauthor{J. V. Zidek}

\address[addr1]{Department of Statistics, University of British Columbia,
2207 Main Mall, Vancouver, BC, Canada V6T 1Z4, \printead{e1}}
\end{aug}

\end{frontmatter}

%% Mainmatter %%

This paper concerns a very topical issue, namely the effect of
preferential sampling the locations at which to measure a spatial
process. The topic was highlighted at and studied by a research group
at the Statistical and Applied Mathematical Sciences Institute (SAMSI)
during its 2009--10 thematic year on spatial statistics, and a number
papers came out of that initiative.

To put this paper in context, some background seems worthwhile.
Selection bias in one form or another has always been an issue in
statistical science, and it has been studied since at least the time
when Horvitz and Thompson proposed their simple but ingenious approach
to unbiasing estimates of finite population averages when sample items
are preferentially selected \citep{horvitz1952}. Survey statisticians
have long since recognized the adverse effect of such bias and the need
to adjust for it when computing their estimates. Biostatisticians have
also been concerned with this issue in the form of response biased
sampling in estimating the relationship between a response $Y$ and a
covariate vector $Z$ when instead of sites human subjects are the units
\citep{scott11}. There, inter-subject dependence is ignored due to its
complexity and the work of \citet{liang1986longitudinal} which allows
that simplification to be made. The responses $Y$ are assumed to be
observed (although that assumption can be relaxed by modeling it) and
subject selection is biased by these responses. In contrast, the
present paper follows \citet{diggle2010geostatistical} and assumes
instead that the role of $Y$ is implicit and seen through the point
process model that ``knows'' $Y=Y(x),$ or rather the latent process
that generates it, through the intensity function $\exp{\{\alpha+
\beta S(x)\}}$, quite a strong assumption. The just cited work in
biostatistics would be of potential relevance in spatial regression, an
important topic in environmental epidemiology, but the effect of
preferential sampling in that domain, especially on the effect on
optimal design as seen in this paper, has not been studied as far as we know.

In geostatistics, spatial dependence can often be of central importance
especially when spatial prediction is of primary interest. The paper by
\citet{diggle2010geostatistical} has awakened interest in a topic that
has been conveniently ignored even by those charged with setting
regulatory standards---where sites may be deliberately sited to detect
the non-compliers and placed where response levels are expected to be
high \citep{guttorp:2010}.

The present paper shares with \citet{diggle2007model}, \citet
{pati2011bayesian} and \citet{gelfand2012effect} the goal of
determining the effect of preferential sampling on statistical
inference, specifically spatial prediction and parameter estimation.
The more recent work of \citet{shaddick2014case} demonstrates fairly
conclusively through a case study that administrators, when left to
their own devices, will select environmental process monitoring sites
preferentially. And the theory in \citet{zidek2014unbiasing} shows
what to do about it, at least if you are an official statistician
reporting annual averages over space of an environmental process
field---that theory relies on the Horvitz--Thompson approach albeit
with estimated selection probabilities. The effect can be quite
dramatic in some years, where the estimated annual averages and numbers
of sites out of compliance with regulatory standards is greatly reduced
relative to their unadjusted counterparts. Finally, in work about to be
submitted, Liu, Shaddick, Zidek, and Cai show that relative risks of
respiratory mortality increase once the effect of preferential sampling
is accounted for, again as seen in a case study. The last three of the
cited works are unlike all the previous papers done in a
spatio-temporal context.

The present paper, done in a spatial context, adds to this series of
papers by showing yet one more impact of preferential sampling, namely
on the all important problem of optimal design. More specifically, the
paper shows where to take additional observations based on spatial data
already collected at preferentially sampled sites, when the
preferential sampling mechanism is known and interest is in a set of
future observations to be collected at the new set of monitoring sites.
This is done with an optimal design theory based on utilities that
measure the benefit of any proposed change to the network, the key
element being the inclusion in the posterior distribution of (the
informative) existing site locations along with the observations at
those sites. Together these two ingredients tell us something about the
relationship between site locations $X$ and responses $Y$.

Although the paper presents a general theory, work is confined to just
utilities based on the predictive variance, in keeping with the
approach to design developed at least 50 years ago by geostatisticians
who added one point at a time to an existing network, placed where
predictive variances are large. However, we would argue that this
approach is too simplistic. An entropy-based approach in the Gaussian
case takes the spatial correlation into consideration so that a site
with high posterior predictive variance may not make it into the new
network because it can be well predicted from some of its newly added
neighbors \citep{le2006statistical}. That is, in fact, what happens in
a real-world application of the entropy-based approach \citep
{ainslie_application_2009} to the redesign of Vancouver's air quality
monitoring network. It would very interesting to see how preferential
sampling would affect new site selection in the rainfall example
considered in this paper when spatial correlation is included. For
recent general discussions of design issues such as these, see \citet
{zidek2010monitoring} and \citet{shaddick2015spatio}.

I have some concern about the accuracy of the approximation given in
Appendix C, especially in the particularly important case where $\beta
\times S$ is not small. It is difficult to intuit how things might
change in that case due to use of the approximation. I wonder if
instead, the Laplace approximation could be used as it has some
theoretical credentials.

The improved results gained by incorporating preferential sampling into
their models relies on knowing the form of model for the preferential
sampling effect. In reality, that model would not be anything like the
right model--selection is complicated by committees using guidelines,
proximity to highways and so on. A starting point in the exploration of
this issue might begin by experimentally adding an $S^2$ term in the
intensity function of the point process.

The work cited above makes the point that covariates need to be
included to see if they can explain the site selection before seeing if
there is a residual role to be played by preferential sampling,
something that is not considered in the theory presented in this paper.

Overall this paper is a welcome addition to the emerging theory of
preferentially sampled spatial monitoring things, and it leaves open
quite a number of new directions to explore. In particular, it offers a
theory for those who must select sites preferentially, for example, to
monitor sources of air pollution.

%% References %%
%

\end{document}